\documentclass[11pt]{article}

\usepackage[a4paper, left=3.3cm,top=3cm,right=3.3cm,bottom=3.7cm]{geometry}
\usepackage{authblk}
\usepackage{framed}
\usepackage{booktabs}
\usepackage{pgfplots}
\usepgfplotslibrary{patchplots}

\usepackage{algorithmic}
\usepackage[ruled]{algorithm}

\usepackage{siunitx}

\usepackage{enumitem}

\newcommand{\trans}{\mathsf{T}}

\newcommand{\argmin}{\mathrm{argmin}}

\usepackage{psfrag}
\usepackage{pifont}

\usepackage{tikz}
\usepackage{tikz-3dplot}

\usepackage{concmath}
\usepackage[T1]{fontenc}

\usepackage{amsthm}
\newtheorem{theorem}{Theorem}
\newtheorem{remark}[theorem]{Remark}

\theoremstyle{definition}

\newcommand{\abs}[1]{\vert#1\vert}

\newcommand\norm[1]{\Vert#1\Vert}

\newcommand{\set}[1]{\{#1\}}
\newcommand{\kl}[1]{(#1)}

\newcommand{\coloneqq}{\triangleq}

\newcommand{\ncoils}{N_c}
\newcommand{\nvoxels}{N_v}
\newcommand{\nsensors}{N_s}
\newcommand{\nact}{m}

\newcommand{\ivoxels}{v}
\newcommand{\icoils}{c}
\newcommand{\iact}{j}
\newcommand{\isensors}{s}

\newcommand{\mnp}{\mathbf{n}}
\newcommand{\znp}{\mathbf{z}}

\newcommand{\Hfield}{\mathbf{H}}
\newcommand{\bfield}{\mathbf{b}}

\newcommand{\Bfield}{\mathbf{B}}

\newcommand{\XX}{\mathbf{X}}
\newcommand{\YY}{\mathbf{Y}}
\newcommand{\y}{\mathbf{y}}

\newcommand{\noise}{\boldsymbol \xi}
\newcommand{\Noise}{\boldsymbol \Xi}

\newcommand{\R}{\mathbb R}

\newcommand{\N}{\mathbb N}

\newcommand{\Ho}{\mathbf H}

\newcommand{\A}{\mathbf A}
\newcommand{\LL}{\mathbf L}

\newcommand{\Io}{\mathcal I}

\newcommand{\Ao}{\mathcal A}
\newcommand{\Lo}{\mathcal L}
\newcommand{\Mo}{\mathcal M}

\newcommand{\sparse}{\boldsymbol \Psi}

\newcommand{\x}{\mathbf x}
\newcommand{\rr}{\mathbf r}
\newcommand{\nv}{\boldsymbol \nu}

\newcommand{\q}{\q}

\newcommand{\Y}{\mathbf Y}

\newcommand{\Hf}{\mathbf H }

\newcommand{\Rcal}{\mathcal R }
\newcommand{\cv}{\mathbf{r}}

\newcommand{\tik}{\mathcal T}
\newcommand{\reg}{\mathcal R}

\usepackage{amsopn}

\DeclareMathOperator{\reshape}{vec}
\DeclareMathOperator{\mat}{mat}

\usepackage{soul}
\usepackage{xcolor}
\colorlet{lred}{red!40}
\colorlet{lgreen}{green!40}
\colorlet{lblue}{blue!40}

\makeatletter
\newcommand*\bigcdot{\mathpalette\bigcdot@{.6}}
\newcommand*\bigcdot@[2]{\mathbin{\vcenter{\hbox{\scalebox{#2}{$\m@th#1\bullet$}}}}}
\makeatother
\newcommand\inner[2]{{#1}\bigcdot {#2}}

\usetikzlibrary{calc,shapes,decorations.markings}

\usepackage{amsmath,amssymb}

\numberwithin{equation}{section}
\numberwithin{figure}{section}

\title{Douglas-Rachford Algorithm for Magnetorelaxometry Imaging using Random and Deterministic Activations}

\allowdisplaybreaks

\author{Markus~Haltmeier}

\affil{Department of Mathematics, University of Innsbruck,\authorcr
Technikerstrasse 13, A-6020 Innsbruck, Austria.\authorcr
{\tt markus.haltmeier@uibk.ac.at}
}

\author{Gerhard~Zangerl}

\affil{Department of Mathematics, University of Innsbruck,\authorcr
Technikerstrasse 13, A-6020 Innsbruck, Austria.\authorcr
{\tt gerhard.zangerl@uibk.ac.at}
}

\author{Peter~Schier}

\affil{Institute of Electrical and Biomedical Engineering,
Private University for Health Sciences, Medical Informatics and Technology (UMIT),
Eduard-Walln\"ofer-Zentrum 1, A-6060 Hall in Tirol, Austria.\authorcr
{\tt peter.schier@umit.at}
}

\author{Daniel~Baumgarten}

\affil{Institute of Electrical and Biomedical Engineering,
Private University, for Health Sciences, Medical Informatics and Technology (UMIT),
Eduard-Walln\"ofer-Zentrum 1, A-6060 Hall in Tirol, Austria.} 
\affil{Institute of Biomedical Engineering and Informatics, Technische Universit\"at Ilmenau,
Gustav-Kirchhoff Str. 2; D-98693 Ilmenau, Germany.\authorcr
{\tt daniel.baumgarten@umit.at}
}

\date{\today}

\begin{document}

\maketitle

\begin{abstract}
 Magnetorelaxometry  imaging is a novel  tool for
 quantitative determination of the spatial distribution  of magnetic nanoparticle
inside an organism. The use of multiple excitation patterns
has been demonstrated  to significantly  improve spatial resolution. However,
increasing the number of excitation patterns is  considerably more time
consuming, because several sequential measurements have  to be performed.

In this paper, we use  compressed sensing in combination with sparse recovery to reduce the total measurement time and to improve spatial  resolution.
For image reconstruction, we propose using the  Douglas-Rachford splitting
algorithm applied to the sparse  Tikhonov functional including a positivity constraint.
Our numerical  experiments  demonstrate that the resulting algorithm
is capable  to accurately recover the  magnetic nanoparticle distribution
from a small number of activation patterns. For example, our algorithm  applied with  10   activations yields half the reconstruction error of  quadratic  Tikhonov regularization   applied with  50 activations, for a tumor-like phantom.


 \end{abstract}

\section{Introduction}

Magnetic nanoparticles (MNP) offer a variety of promising biomedical applications.
For example, they can be used as agents for drug delivery or hyperthermia,
where the aim is to heat up specific regions inside a biological specimen~\cite{Wiekhorst2012}.
These applications require quantitative knowledge of the  magnetic nanoparticles
 distribution for safety and efficacy monitoring.
 In this paper, we consider  magnetorelaxometry  (MRX) imaging, which is a novel and promising non-invasive technique to spatially resolve the location  and concentration of MNPs in
 vivo~\cite{Baumgarten2008,Liebl2014}. It beneficially combines a highly sensitive magnetic measurement technology for MNP imaging with a broad range of parameters and has the potential to image particle distributions in a comparably large volume~\cite{coene2017multi}.

In MRX, the magnetic moments of the magnetic nanoparticles are aligned by a magnetic field generated by excitation coils~\cite{Liebl2015}. Therewith, a magnetization can be measured from the MNP. After switching off the excitation field, the decay of magnetization (relaxation) is recorded, typically by
SQUIDs (superconducting quantum interference device), yielding information about the particle concentration and related properties. If the relaxation is measured by a sensor array, the particle distribution can be reconstructed by inverse imaging methods~\cite{Baumgarten2010, Liebl2014}.

In multiple excitation MRX (ME-MRX), several inhomogeneous activation fields generated by an array of excitation coils are applied. With such a coil array, a variety of inhomogeneous excitation fields can be generated, for example, by switching on the coils in a sequential manner~\cite{Steinhoff2010}.
First experimental realizations of ME-MRX with sequential activation and least squares estimation for imaging
have been obtained in~\cite{Liebl2014}.
ME-MRX using sequential coil activation and standard image reconstruction techniques requires a large number of measurement cycles and thus a considerable time for data acquisition. Recently, different approaches have been proposed using advanced excitation schemes~\cite{Coene2012, crevecoeur2012advancements, Baumgarten2015}.

A strategy to maintain high resolution while reducing measurement time
is via compressed sensing (CS), a  new sensing  paradigm \cite{CanRomTao06b,Donoho2006compressed} that allows to capture a high-resolution image (or signal) by using fewer measurements than predicted by Shannon's sampling theorem. CS replaces point-measurements by general linear  measurements,
where each measurement consists of a linear combination of the entries
of the image  of interest. Recovering the original image is highly under-determined.
The standard theory of CS predicts  that under suitable assumptions on the image
(sparsity) and the measurement matrix (incoherence), stable image
 reconstruction is  nevertheless possible. CS has led to several new sampling strategies in medical imaging,
 for example, for speeding up MRI data acquisition~\cite{lustig2008compressed},
accelerating photoacoustic tomography~\cite{haltmeier2016compressed},
or completing under-sampled CT images~\cite{chen2008prior}.  While these applications involve
relatively  mildly ill-conditioned problems, ME-MRX constitutes a severely ill-conditioned  or ill-posed problem
\cite{baumgarten2013magnetic,focke2018inverse}.
  No standard approaches for image  reconstruction
combining compressive measurements and such severely ill-posed problems exist.

In this work, we investigate CS techniques for accelerating  ME-MRX. For that purpose, we consider random activations of the coils as well as sparse  deterministic activation schemes.
In order to image the MNP distribution from the CS data we
propose a sparse reconstruction framework. We use $\ell^1$-Tikhonov regularization together with a positivity constraint on the set of reconstructed
 MNP distributions, and use the Douglas-Rachford splitting algorithm \cite{combettes2011proximal} for its minimization.
 We evaluate the performance of the resulting sparse recovery algorithm  using different phantoms in dependence of the number of activation patterns. In any of the cases,
 the sparse recovery algorithm reduces the root mean square error (RSME)
by more than $\SI{50}{\percent}$. For two of the three phantoms (CS and smiley),
the deterministic scheme  outperforms the random pattern   by around $\SI{10}{\percent}$ in terms of  reduced RSME. For the third phantom (tumor), all tested sampling patterns
roughly  give the same results. We therefore can propose the  sparse sampling pattern
among the tested schemes, in combination with the proposed sparse recovery framework.

\section{The forward model in ME-MRX imaging}
\label{sec:me-mrx}

This section gives a precise formulation  of the forward problem
of  ME-MRX and CS.  Throughout this text we work with a  discrete model~\cite{Wiekhorst2006,Baumgarten2008,Richter2010}.
A continuous model has recently  been presented in
\cite{focke2018inverse}.

\begin{figure}[htb!]
	\center
	\begin{tikzpicture}[scale=1.2,
trans/.style={thick,<-,shorten >=2pt,shorten <=2pt,>=stealth}
]	
	\draw[step=0.25cm,black,line width = 1.5pt] (0,0) grid (3,3);
	\foreach \x in {0.125,  0.375,...,2.875}
	\foreach \y in {0.125,  0.375,...,2.875}
	{
		\draw[fill = black ] (\x, \y) circle (1pt);
	}

	\foreach \y in {-0.5,0, 0.5,...,3}
	{
		\draw[trans] (3.6,\y) -- (4,\y) node[midway,right] {};
		\filldraw[even odd rule,inner color=red,outer color=red, opacity = 0.6]
		(4,\y) circle (0.0)
		(4,\y) circle (0.1);

		\draw[trans] (-0.6,\y) -- (-1,\y) node[midway,right] {};
		\filldraw[even odd rule,inner color=red,outer color=red, opacity = 0.6]
        (-1,\y) circle (0.0)
		(-1,\y) circle (0.1);

	}

	\foreach \x in {-0.5,0,0.5,1,1.5,2,2.5,3,3.5}
	{
		\draw[trans] (\x, -0.6) -- (\x, -1) node[midway,right] {};
		\filldraw[even odd rule,inner color=red,outer color=red, opacity = 0.6]
		(\x, -1) circle (0.0)
		(\x, -1) circle (0.1);

	}

		\filldraw[even odd rule,inner color=red,outer color=red, opacity = 0.6]
 	    (-1, -1) circle (0.0)
		(-1, -1) circle (0.1);
	
		\filldraw[even odd rule,inner color=red,outer color=red, opacity = 0.6]
 	    (4, -1) circle (0.0)
		(4, -1) circle (0.1);

	\draw [step = 0.25, draw = black,  fill= cyan, opacity = 0.6]  (-2,3.5) grid (5,4); \draw[fill = cyan, opacity = 0.6] (-2,3.5) rectangle (5,4);  	
	\node at (4.6, 1) {\rotatebox{-90}{  activation coils   }};
	\node at (1.5, 4.3) {\rotatebox{0}{   sensors  }      };
	\node at(1.5, -0.3) {   voxels    };
	\end{tikzpicture}
\caption{Illustration  of a simplified two-dimensional ME-MRX setup consisting of an array of activation
coils and two layers of sensors arranged around the volume of interest.
The coil locations $\rr_\isensors$ are indicated by red circles, and the corresponding normal vectors $\nv_\isensors$ by black arrows.}\label{fig:MRX-setup}
\end{figure}
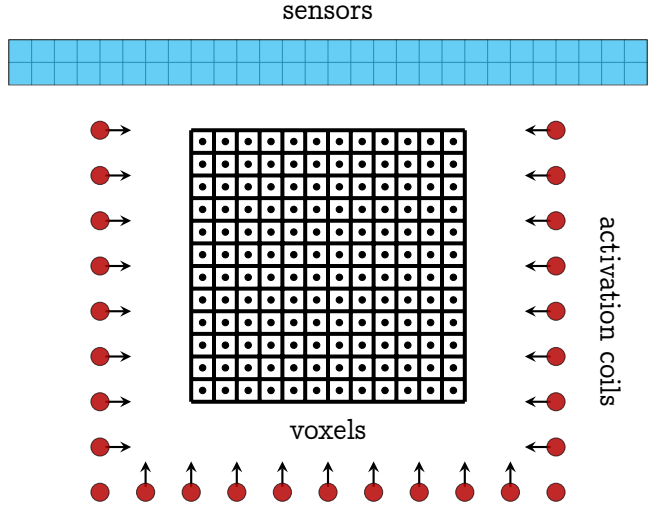

\subsection{Signal generation}

In  MRX, an activation field is  applied to a region of interest containing MNPs.
We assume  this region to be divided into
a number of $\nvoxels$ quadratic voxels, each represented
by its midpoint located at $\rr_\ivoxels$  and containing a concentration $n(\rr_\ivoxels)$
of MNPs  (compare Figure~\ref{fig:MRX-setup}).
Due to the presence of the applied field, the  MNPs align up and,
after the activation field is switched off, produce a relaxation signal.
According to the Biot-Savart law, the overall contribution of the MNP concentrations
in the individual voxels to the measured signal at sensor location $\rr_\isensors$
is given  by~\cite{Liebl2014}
\begin{equation}\label{eq:forw_sol}
b^{\Ho} (\rr_\isensors)  = \frac{\mu_{0}}{4 \pi} \; \nv_\isensors \bigcdot  \sum_{\ivoxels=1}^{\nvoxels}
\biggl( \frac{3  \rr_{\isensors,\ivoxels} \otimes \rr_{\isensors,\ivoxels}  }{| \rr_{\isensors,\ivoxels}|^5}
 - \frac{1}{|\rr_{\isensors,\ivoxels}|^3} \biggr) \;
\Ho (\rr_\ivoxels) \;  n(\rr_\ivoxels)  \quad \text{ for } \isensors \in \set{1, \dots, \nsensors}\,.
\end{equation}
Here $\Ho$ is  the  activation field, $\rr_{\isensors,\ivoxels} = \rr_\isensors - \rr_\ivoxels$
is the vector joining $\rr_\isensors$ and $\rr_{\ivoxels}$,  $\nv_\isensors$
is the normal vector of the sensor at location $\rr_\isensors$, and $\otimes$ and $\bigcdot$ are  used to
denote the tensor and inner product, respectively. Recall that the tensor product of two column vectors
$\mathbf{a}, \mathbf{b} \in \R^3$ is defined  by $\mathbf{a} \otimes \mathbf{b} \coloneqq \mathbf{a} \mathbf{b}^\trans \in \R^{3 \times 3}$
which results in a matrix, whereas the inner product $\mathbf{a}  \bigcdot \mathbf{b} \coloneqq \mathbf{a}^\trans \mathbf{b} \in \R$
results in a scalar.

Assuming the concentrations $n(\rr_\ivoxels)$ and the activation field  $\Ho$ to be known,
the  measured data can be computed by~\eqref{eq:forw_sol}.
Collecting  the individual measurements in a vector $\bfield^{\Ho}   = (b^{\Ho}  (\rr_1) , \dots, b^{\Ho}  (\rr_{\nsensors}))^\trans$,
we obtain the linear equation
\begin{equation}\label{eq:leadfield}
	\bfield^{\Ho}  = \Lo^{\Ho} \, \mnp  \in \R^{\nsensors}\,.
\end{equation}
Here  $ \mnp = (n(\rr_1), \dots, n(\rr_{\nvoxels}))^\trans \in \R^{\nvoxels}$ is the vector of
MNP concentrations,  and   $ \Lo \in \R^{\nsensors \times \nvoxels}$ is the system matrix
having entries
\begin{equation}\label{eq:lead}
\Lo^{\Ho}_{\isensors,\ivoxels}  \triangleq \frac{\mu_{0}}{4 \pi} \; \inner{\nv_\isensors}{
\biggl( \frac{3  \rr_{\isensors,\ivoxels} \otimes \rr_{\isensors,\ivoxels}  }{| \rr_{\isensors,\ivoxels}|^5}
- \frac{1}{|\rr_{\isensors,\ivoxels}|^3} \biggr)
\Ho (\rr_\ivoxels) }
\end{equation}
as derived from relation \eqref{eq:forw_sol}. The matrix $\Lo^{\Ho}$ is called  Lead field matrix corresponding to the excitation field  $\Ho$.
Equations~\eqref{eq:leadfield}, \eqref{eq:lead}  constitute the standard discrete
forward model  of  MRX  using a single activation field.

\begin{figure}
\centering \includegraphics[width=0.7\textwidth]{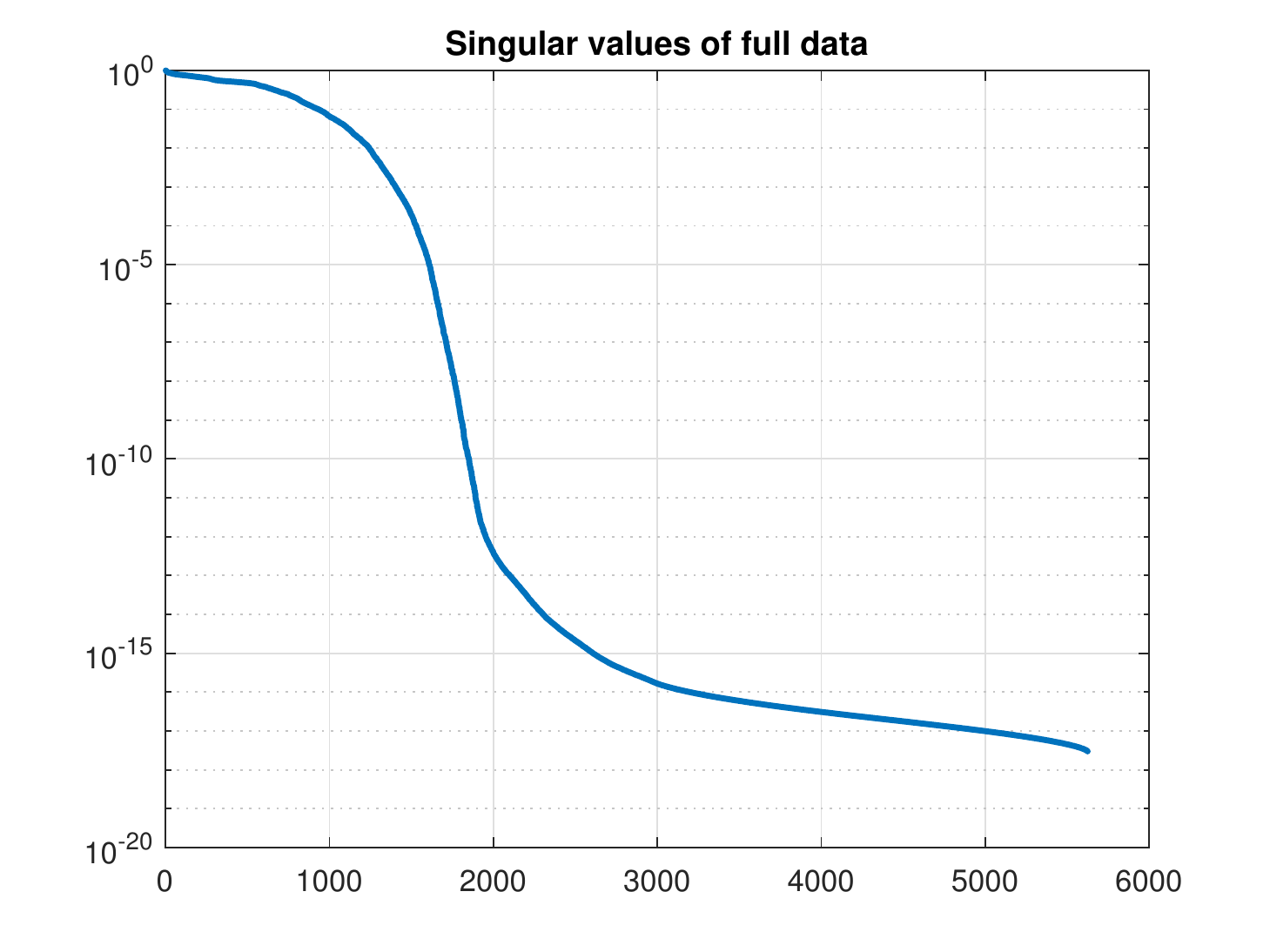}
\caption{\textbf{Singular values  of  the Lead field matrix on a logarithmic scale.} \label{fig:singular}
The  Lead field matrix  uses $110 \times  120 = 13200$
measurements for $75^2 = 5625$  unknowns.
The  rapid decay of the singular values is due to the severe
ill-conditioning of the matrix. Increasing the number
of measurements results in a similar rapid decay
which indicates the inherent ill-posedness of the inverse problem.}
\end{figure}

\subsection{Sequential coil activation}

In ME-MRX, the volume is (partially) surrounded by an array of excitation coils and sensor arrays (see Figure~\ref{fig:MRX-setup}). The coils can be steered  to generate different types of activation fields.
In the following we denote by $\Hfield_1, \dots, \Hfield_{\ncoils}$  the magnetic
fields induced by individual activation of the coils. Further, we write  $\bfield_\icoils \triangleq \bfield^{\Hfield_{\icoils}} \in \R^{\nsensors}$ and
$\Lo_{\icoils} \triangleq \Lo^{\Hfield_{\icoils}} \in \R^{\nsensors\times \nvoxels }$  for the data
and the Lead field matrix according to \eqref{eq:leadfield} and  \eqref{eq:lead},
corresponding to the  activation field  $\Hfield = \Hfield_\icoils$.
The  measurements from all activations
can be combined to a single equation of the form
\begin{equation}\label{eq:Mlied}
\bfield = \Lo \, \mnp \in   \R^{\ncoils  \nsensors}  \,,
\end{equation}
where  
\begin{align}\label{eq:Mlied2}
\bfield \triangleq
\begin{bmatrix}
\bfield_1 \\
\rule{1em}{0.5pt}\\
\vdots \\
\rule{1.5em}{0.5pt}\\
\bfield_{\ncoils}
\end{bmatrix} \in   \R^{\ncoils  \nsensors}
\quad \text{ and } \quad
\; \Lo  \triangleq
\begin{bmatrix}
\Lo_1 \\
\rule{1em}{0.5pt}\\
\vdots \\
\rule{1em}{0.5pt}\\
\Lo_{\ncoils}
\end{bmatrix}  \in   \R^{\ncoils  \nsensors \times \nvoxels } \,.
\end{align}
We will refer to $\bfield$ in \eqref{eq:Mlied} as full activation data.
Evaluating   \eqref{eq:Mlied}  constitutes the forward model of ME-MRX.
The corresponding inverse problem consists in determining the MNP distribution $\mnp \in \R^{\nvoxels}$ from the data vector $\bfield$
that is additionally  corrupted  by noise.
Note that  Eq.~\eqref{eq:Mlied}  is known to be severely ill-conditioned as its singular
values are rapidly decreasing  (see \cite{baumgarten2013magnetic};
 an analysis in  the infinite-dimensional setting has been performed in \cite{focke2018inverse}).
In Figure {\ref{fig:singular}} we plot the
singular values of the forward matrix for full data.
The singular values decay rapidly. This reflects  the severe ill-conditioning
 and implies that only a small number of singular components can be robustly
 reconstructed.

The  multiple coil setup has been realized in \cite{Baumgarten2015,Steinhoff2010}
and shown to  significantly
improve the  spatial resolution compared to a single coil activation~\cite{Steinhoff2010}.
However, consecutive  activations lead to a more time consuming measurement
process.  In order  to accelerate  data acquisition and to improve imaging resolution, we use CS techniques   and
develop  a Douglas-Rachford based sparse reconstruction algorithm.

\subsection{Compressive coil activations}

The basic idea to employ  CS for ME-MRX is to use $\nact$ specific
coil  activations, with $\nact  \ll \ncoils$, instead of activating all
$\ncoils$ coils in a sequential manner.    Because of
linearity,  the data corresponding  to the $ \iact$th
random  activation  pattern are given by
\begin{equation}\label{eq:csmeas}
	\y_\iact
	=  \sum_{\icoils = 1}^{\ncoils}    a_{\iact, \icoils} \bfield_\icoils
	=  \sum_{\icoils = 1}^{\ncoils}    a_{\iact, \icoils} ( \Lo \mnp)_\icoils + \noise_\iact
	\quad  \text{ for } \iact = 1, \dots, \nact   \,.
\end{equation}
Here $\noise_\iact$ models the error in the data and   $a_{\iact, \icoils}$ is the
contribution of coil $\ncoils$ to the   $\iact$th activation pattern.
The measurement matrix
\begin{equation} \label{eq:matrix}
\A  \coloneqq
\begin{pmatrix}
a_{1, 1}& \cdots& a_{1, \ncoils} \\
\vdots & & \vdots \\
a_{\nact, 1}& \dots& a_{\nact, \ncoils}
\end{pmatrix}  \in \R^{\nact \times \ncoils}
\end{equation}
represents all activation patters.
Typical choices for $\A$ are
Bernoulli or Gaussian  matrices, since they are known to guarantee stable
recovery of sparse signals \cite{foucart2013mathematical}.
Such kind of measurements can be  naturally realized for ME-MRX,
by simultaneous activation of  several coils.
The corresponding image reconstruction   problem  consists in
recovering the   MNP distribution $\mnp$ from the  data in
\eqref{eq:csmeas}.

\subsection{Matrix formulation of the reconstruction problem}

In order to write the reconstruction  problem  in a compact form we introduce some additional
notation.   First, we define the CS measurement vector   as
\begin{equation}
\y \triangleq
\begin{bmatrix}
\y_1 \\ \rule{1.5em}{0.5pt}\\ \vdots \\ \rule{1.5em}{0.5pt}\\ \y_{\nact}
\end{bmatrix} \in \R^{\nact \nsensors} \,.
\end{equation}
Second,   we introduce the  vectorization or
reshaping operator $\reshape_{\ncoils}$, that   takes a matrix  to a  vector  whose
block entries are equal to the transposes of   the  rows of the matrix:
\begin{equation}
\reshape_{\ncoils} \colon \R^{\ncoils \times  \nsensors}  \to  \R^{ \ncoils  \nsensors }
\colon  \XX  \mapsto
\begin{bmatrix}
\XX_{1, -}^\trans
\\
\rule{1.5em}{0.5pt}\\
\vdots
\\
\rule{1.5em}{0.5pt}\\
\XX_{\ncoils, -}^\trans
\end{bmatrix} \,.
\end{equation}
We denote by $\mat_{\ncoils}$ the inverse reshaping
operation that maps  a  vector in $\R^{ \ncoils  \nsensors }$  to a matrix in $\R^{ \ncoils \times
 \nsensors }$.
Further,  we write $\YY \coloneqq  \mat_\nact \y$
and $\Bfield \coloneqq  \mat_{\ncoils} \bfield$.

With the above notations, the CS data  in \eqref{eq:csmeas}
can be written in the compact form
\begin{equation} \label{eq:singleCS2} 
\Y= \A \,  \LL \, \mnp  + \Noise   \in \R^{   \nact  \times \nsensors}   \,,
 \end{equation}
where $\LL$ is the reshaped Lead field matrix  defined by the  property
$\LL \mnp  = \mat_{\ncoils}  (  \Lo \mnp)$ and $\Noise $ denotes the noise matrix.
With the Kronecker (or tensor)  product  $ \Ao \coloneqq \A  \otimes \Io \in \R^{(\nact \nsensors) \times (\ncoils \nsensors)}$
between the CS activation matrix  $\A$ and the identity matrix $\Io  \in \R^{\nsensors \times \nsensors} $,
equation~\eqref{eq:singleCS2} can further be rewritten in the form
\begin{equation} \label{eq:singleCS3} 
\y  =  \Ao \,   \Lo \, \mnp + \noise \in \R^{ \nact  \nsensors }  \,.
\end{equation}
The image reconstruction task in CS ME-MRX consists in recovering
 the MNP distribution $\mnp$ from equation~\eqref{eq:singleCS2} or, equivalently, from
equation~\eqref{eq:singleCS3}.

 \section{Douglas-Rachford algorithm for CS MRX imaging}
 \label{sec:rec}

In order to recover  the MNP distribution,
the inverse problem  {\eqref{eq:singleCS2}}  (or {\eqref{eq:singleCS3}})
has to be solved which is known to be severely ill-posed,
see Figure~{\ref{fig:singular}}.
We note that no standard approach  for  image
 reconstruction combining
compressive measurements and severely ill-posed problems exist.
This section  is also of interest from a general  perspective on
CS for inverse  problems.

\subsection{Sparse Tikhonov regularization}

We consider \eqref{eq:singleCS3} as a single  inverse
problem with system matrix $ \Mo = \Ao\Lo $.   To address instability and non-uniqueness and to
incorporate  prior knowledge, we use a sparse regularization approach. For that purpose, we
minimize the generalized Tikhonov  functional~\cite{scherzer2009}
 \begin{equation} \label{eq:tik1}
	\tik_{ \alpha, \beta }^{\Mo} (\mnp)  \triangleq
	\frac{1}{2} \norm{ \y  - \Mo \mnp }_2^2
	+ \alpha  \norm{ \sparse  \mnp}_1
	+ \beta   \reg(\mnp)  \,.
\end{equation}
Here $\sparse \colon \R^{\nvoxels} \to \R^{d}$ is a transform that sparsifies the
MNP concentration (with the $d \in \N $ the dimension of the transformed domain),
$\Rcal \colon \R^{\nvoxels} \to [0, \infty)$  is an additional
 regularizer that incorporates additional prior knowledge  about the MNP
 distribution (such as positivity and other convex constraints) and
$\norm{ \, \cdot  \, }_p$ is the standard $\ell^p$-norm defined by
\begin{equation*}
\norm{ \x  }_p \triangleq \sqrt[p]{ \sum_{i=1}^d \abs{x_i}^p } \quad \text{ for }
\x =(x_1, \dots, x_d)\in \R^d \,.
\end{equation*}
The non-negative parameters $\alpha$ and $\beta$  allow to
balance between the data consistency term $\frac{1}{2} \norm{ \y  - \Mo \mnp }_2^2$,
the  sparsity prior $ \norm{ \sparse  \mnp}_1$  and  the additional regularizer
$\reg(\mnp)$.

We call $\sparse$ a sparsifying transform if $\sparse \mnp$ can well be approximated by $k$-sparse vectors $\x$, defined by the property that $\set{i \mid  x_i \neq 0 }$ has at most $k$ elements. This approximation be can be quantified using
the $k$-term approximation error,
\begin{equation} \label{eq:sterm}
 	\sigma_{k} (\mnp)
	\coloneqq \inf \set{ \norm{\mnp - \mnp_k}_1
	\mid  \mnp_k \textnormal{ is $k$-sparse} }.
\end{equation}
The $k$-term approximation error  $\sigma_{k} (\mnp)$
appears as additive term  in standard error estimates in
compressed sensing theory. For a mathematically precise
discussion of  the   $k$-term approximation error we refer
to \cite{foucart2013mathematical}.

\begin{remark}[Recovery theory for \eqref{eq:tik1}]\label{rem:recovery}
The $k$-restricted isometry property ($k$-RIP)
of $\Mo = \Ao\Lo$ (after appropriate scaling) is defined as the
smallest number $\delta_k$ such that for all $k$-sparse
vectors $\sparse \mnp$ we have
\begin{equation} \label{eq:RIP}
(1-\delta_k) \norm{ \mnp }_2^2 \leq \norm{ \Mo \mnp }_2^2\leq(1+\delta_k)  \norm{ \mnp }_2^2   \,.
\end{equation}
Roughly spoken, CS theory \cite{CanRomTao06b,foucart2013mathematical}
predicts uniform stable recovery with  \eqref{eq:tik1}
for $\beta =0$ and in the sense that  it stably recovers any $s$-sparse vector provided that the $s$-RIP constant of $\Mo$ is sufficiently small.

Due to the severe ill-conditioning of the Lead field matrix $\Lo$, the $s$-RIP  is   expected to be at most satisfied when the number $\nvoxels$ of voxels is small.
Figure~{\ref{fig:singular}} shows that the singular values of $\Lo$
with  $\nvoxels = 5625$ are rapidly decaying and therefore the estimate
$(1-\delta_k) \norm{ \mnp }_2^2 \leq \norm{ \Mo \mnp }_2^2$ for reasonable
$\delta_k$  can only be satisfied for signals mainly formed by the first singular vectors.
This implies that even in the case of full measurement data, where $\Ao$ is the identity matrix, uniform stable reconstruction
 of all sparse vectors is  impossible. To obtain quantitative error estimates for
\eqref{eq:tik1},  stable recovery results for individual elements  \cite{Fuc05b,GraHalSch11,haltmeier2013stable} are a  promising  alternative.
Such theoretical investigations are beyond the scope of this paper and an interesting line of future research.

Standard CS theory requires orthogonality of the sparsifying transform, which  means  that  $\sparse^\trans \sparse$ equals the identity matrix.
Total variation (TV)-regularization \cite{SchGraGroHalLen09}
 uses the gradient as sparsifying transform which is non-orthogonal.
 As in other compressed sensing imaging applications~\cite{krahmer2014stable,poon2015role}
 we observed TV  to outperform orthogonal  sparsifying
 transforms. We will therefore focus on TV and, due to space limitations,
 no results for orthogonal bases like wavelets are included.
\end{remark}

\subsection{Douglas-Rachford algorithm}

In order to minimize \eqref{eq:tik1}, we propose using
the  Douglas-Rachford minimization algorithm, which is a
backward-backward type splitting method for minimizing  the sum
$F+G$ of two  functionals $F$ and $G$. For our  purpose  we take
\begin{align}
F(\mnp)  & \triangleq \frac{1}{2} \norm{ \y  - \Mo \mnp }_2^2 \\
G(\mnp)  & \triangleq  \alpha  \norm{ \sparse \mnp}_1 + \beta \reg(\mnp) \,.
\end{align}
Minimizing \eqref{eq:tik1} by the  Douglas-Rachford algorithm, generates a sequence $(\mnp_k)_{k \in \N}$ of estimated MNP distributions
 and auxiliary sequences $(\znp_k)_{k \in \N}$, $(\tilde\znp_k)_{k \in \N}$
as described in Algorithm~\ref{alg:dr}. See Section~\ref{sec:fwdnum}, for a discussion of the role of the parameter
$\mu$.

\begin{algorithm}
\caption{Proposed \label{alg:dr} sparse reconstruction algorithm for CS ME-MRX.}
\begin{algorithmic}[1]
\STATE Select  $s \in (0,2)$, and $\mu >0$
\STATE Initialize $\znp_0 = 0$
\FOR{ $k=1, \dots, N_{\rm iter}$}
\STATE $\mnp_k \gets (\Mo^\trans \Mo + \mu  \Io)^{-1} (\Mo^\trans \y + \mu \znp_k )$
\STATE $ \tilde \znp_k\gets \argmin_{ \tilde \znp}  \frac{\mu}{2}\norm{(2 \mnp_k - \znp_k) - \tilde\znp}_2^2 + G(\tilde\znp)$
\STATE $ \znp_{k+1} \gets \znp_k
+  s \kl{  \tilde \znp_k  - \mnp_k}$
\ENDFOR
\end{algorithmic}
\end{algorithm}

The Douglas-Rachford  algorithm   is a splitting type   algorithm for minimizing $F+G$ which alternately performs updates according  for $F$ (in our case, the  residual functional or data fitting term)
and  $G$ (in our case, the regularizer). Because the regularizing  term is non-smooth, standard methods such as Newton type methods  cannot be applied.
Splitting type methods are a natural choice in this case, and the implicit
step in line 5 in Algorithm~{\ref{alg:dr}} accounts for the non-smoothness of  $G$.
Other splitting algorithms that are well suited to treat the  non-smooth regularizer are
the forward-backward splitting algorithm \cite{combettes2011proximal}
and the Chambolle-Pock algorithm \cite{ChaPoc11}.
These algorithms use explicit gradient updates for the data fitting term
requiring a matrix inversion. Therefore, they are not applicable to ME-MRX imaging.
The Douglas-Rachford algorithm performs an implicit update  for the data fitting term
(line 4 Algorithm~{\ref{alg:dr}}), which is  recognized as quadratic Tikhonov regularization
step for the given equation and potentially accelerates the iteration.

For our  numerical experiments, we  take $\sparse = \nabla $ as the discrete
gradient operator as appropriate sparsifying transform for piecewise smooth MNP distributions.
The additional regularizer  is taken as the indicator function  $\reg = I_C $
of the convex set
\begin{equation*}
C  \triangleq \set{\mnp \mid  \forall \rr_\ivoxels \colon 0  \leq n(\rr_\ivoxels) \leq n_{\rm max} }
\quad \text{ for some bound }  n_{\rm max} \,,
\end{equation*}
defined by $I_C(\mnp) = 0 $ if  $\mnp \in C$ and $I_C(\mnp) =\infty $. It  guarantees
non-negativity and boundedness.    With the above choices, the Tikhonov  functional
\eqref{eq:tik1} reduces to total variation (TV) regularization with positivity constraint.
 The minimization of
$\frac{\mu}{2}\norm{\mnp - \znp}_2^2 + \norm{ \nabla \znp}_1
+ i_C(\znp) $  required for implementing Algorithm~\ref{alg:dr} is a TV denoising step and
again performed  by the Douglas-Rachford algorithm using the decomposition in $F(\znp) = \frac{\mu}{2}\norm{\mnp - \znp}_2^2 + \norm{ \nabla \znp}_1$ and $G(\znp) = i_C(\znp)$.

Under the reasonable assumptions that the regularizer  $G$ is  lower semicontinuous and convex, and that the sparse Tikhonov functional $\tik_{\alpha, \beta}^{\Mo}$ is coercive, the sequence   $(\mnp_k)_{k\in \N}$
generated by Algorithm~\ref{alg:dr} is known to converge to a minimizer of
\eqref{eq:tik1}; see~\cite{combettes2011proximal,svaiter2011weak}.
Recall that the functional $G$ is called lower semicontinuous  if
$G(\mnp) \leq \lim_{k \to \infty} G(\mnp_k)$ for any  $\mnp$ and any
sequence $(\mnp_k)_{k\in \N}$ converging to $\mnp$.
Note that Algorithm~\ref{alg:dr}  performs implicit steps  with respect to the  residual functional $F = \frac{1}{2} \norm{\Mo \mnp -\y}_2^2$, which  we found to have much faster convergence than forward-backward splitting algorithm \cite{combettes2011proximal}
\begin{align}
 \znp_k  &  \triangleq  \mnp_k -  \Mo^\trans (\Mo \mnp_k -\y )  \\
 \\mnp_k  &  \triangleq \argmin_{\mnp}
 \frac{\mu}{2}\norm{\znp_k - \mnp}_2^2 + G(\mnp)\,
\end{align}
where $\znp_k$ are auxiliary quantities and  $\mnp_k$ the reconstructed MNPs.

\begin{figure}[htb!]
	\centering
	\includegraphics[width = 0.7\columnwidth]{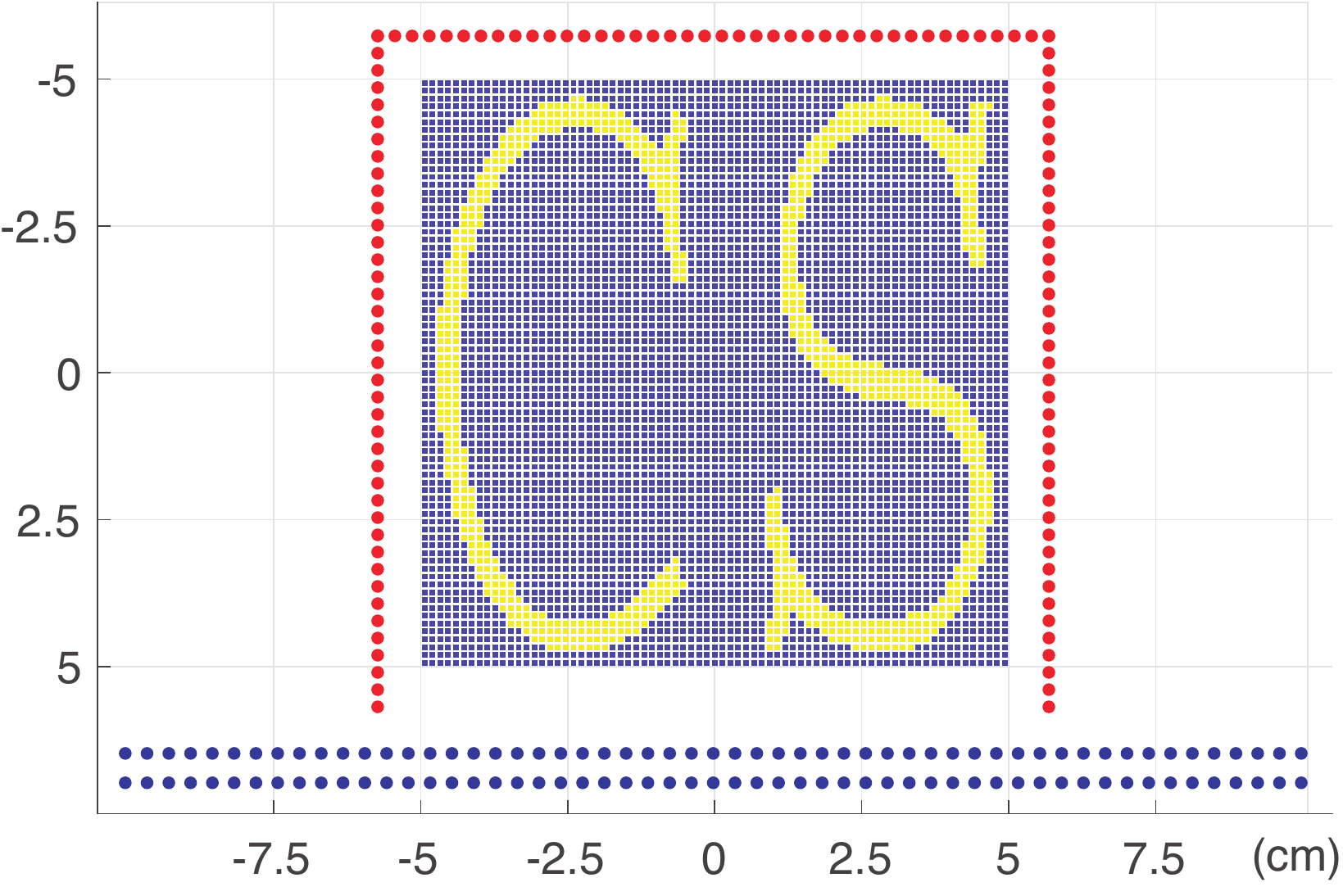}
	\caption{\textbf{Full data setup for the numerical simulations.}
	The phantom is contained in a
	square region of interest, the sensors are arranged in two parallel layers (blue dots),
	and the coils for full activation are  arranged around the phantom in
	 U-form (red dots).}
	\label{fig:NUMsetup}
\end{figure}

\section{Numerical results}\label{sec:num}
\label{sec:num}

For the following numerical simulations, we use a data setup similar to the realization in \cite{Baumgarten2015}. For simplification, we consider a two-dimensional setup representing one voxel plane, containing two parallel arrays of  detector elements
(one measuring the horizontal $(1,0)$ and one  measuring the vertical $(0,1)$ component)
located outside a quadratic region of interest. Circular shaped activation coils are arranged in U-form around the region of interest, all having the same normal vector $(0,1)$.
The used arrangement of  sensor and coil locations for the full measurement data  setup  is shown in  Figure~\ref{fig:NUMsetup}.
For our simulations, we choose a discretization of imaging space into $\nvoxels = 75^2$
voxels covering a region of interest of $[-5,5] \times [-5,5]$  \si{cm^2}.   The   data are generated  for  $\nsensors = 110$
positions and full measurement data  correspond to $\ncoils = 120$ activation coils outside the region of interest.

\begin{figure}[htb!]
	\begin{center}
	\includegraphics[width = 0.9\columnwidth]{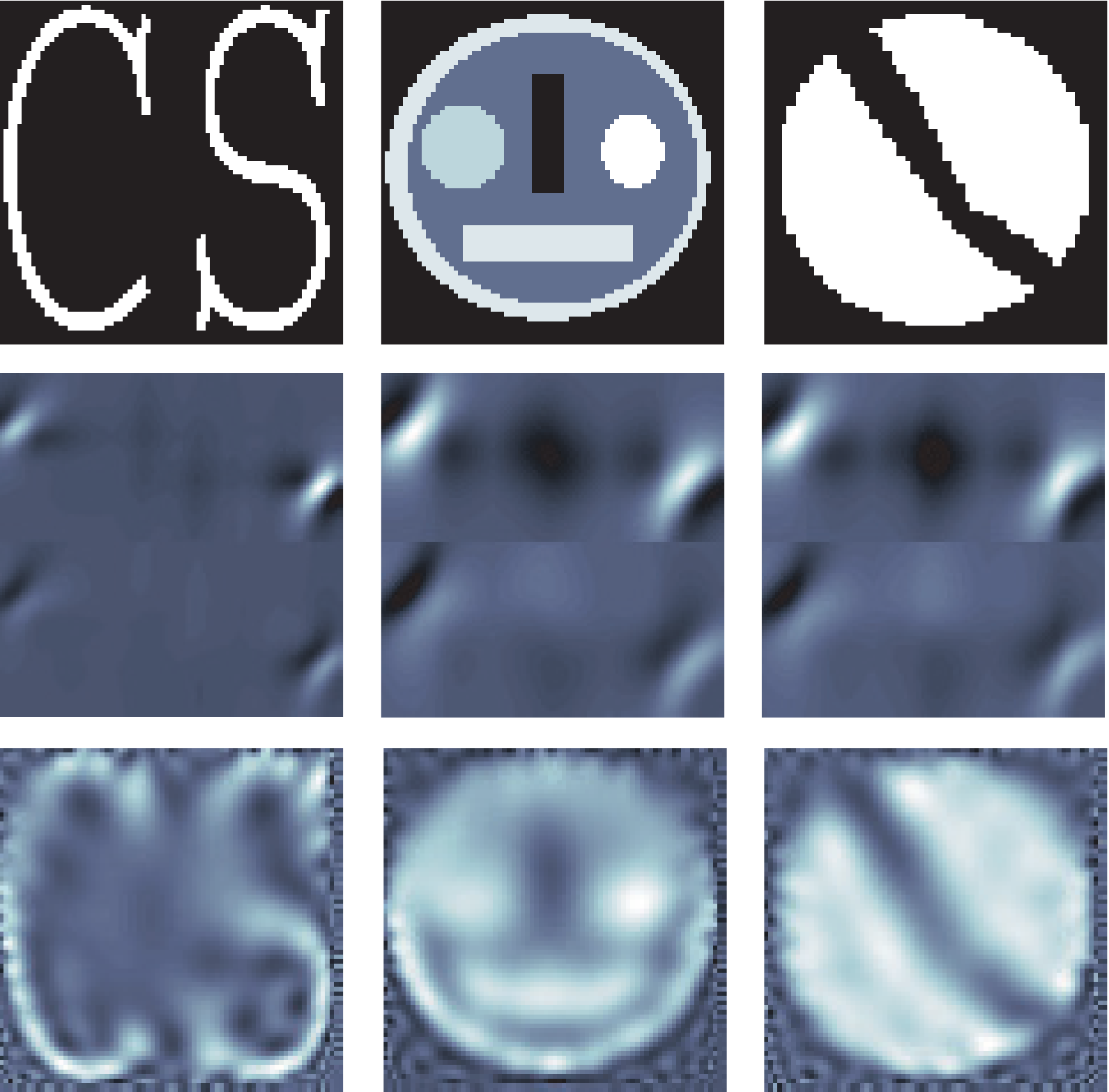}
	\caption{\textbf{Phantoms, data and full measurement data reconstruction.}
	Top row:
	CS-phantom (left), smiley-phantom (middle) and tumor-phantom
	mimicking cancerous tissue with included blood vessel (right).
	Middle row: Corresponding  full measurement data.
	Bottom row: Reconstructions from full measurement data using standard quadratic
	Tikhonov regularization.\label{fig:phantoms}}
	\end{center}
\end{figure}

\subsection{Forward computations}
\label{sec:fwdnum}

To set up the forward model   \eqref{eq:leadfield}, \eqref{eq:Mlied}, \eqref{eq:Mlied2},
we first have to compute the  magnetic fields generated by  each coil.
For that purpose, we follow the approach of \cite{Baumgarten2010,Baumgarten2015}, where any circular activation coil is approximated by short line segments, each carrying  the same current $I_0$. Every coil is modeled by 45 line segments for numerical calculations of the magnetic field. This is considered to be a sufficiently accurate approximation since the deviation of only 36 line segments from the analytical solution was shown to be already below  \SI{1}{\percent} in \cite{liebl2011bild}.
Using this approximation, the  induced magnetic field at the voxel center
 $\rr_\ivoxels$ can be computed by (see~\cite{hanson2002compact})
\begin{equation}
\Hf (\rr_\ivoxels) \simeq  \frac{1}{4\pi} \sum_{i = 1}^{45}    \frac{\abs{ \cv_{1,\ivoxels,i} }
+ \abs{ \cv_{2,\ivoxels,i} } }{ \abs{\cv_{1,\ivoxels,i} } \,  \abs{ \cv_{2,\ivoxels,i} }}
\cdot \frac{
\cv_{1,\ivoxels,i}  \times \cv_{2,\ivoxels,i}   }{ \abs{ \cv_{1,\ivoxels,i} } \, \abs{\cv_{2,\ivoxels,i}}  +  \cv_{1,\ivoxels,i} \bigcdot \cv_{2,\ivoxels,i}} \; I_0  \,,
\end{equation}
where    $\cv_{1,\ivoxels,i}$ and $\cv_{2,\ivoxels,i}$  are  the
 distance vectors between the voxel center $\rr_\ivoxels$  and the
 beginning and end points of the $i$th line segment, respectively. For the
 presented numerical computations, we use a coil
 diameter  of $\SI{1}{\mu m}$, illustrating an almost point like coil.
 Having computed the activation field $\Hf (\rr_\ivoxels)$,  we compute
the entries of the Lead field matrix according to  \eqref{eq:lead}.
By activating the coils sequentially, we obtain the full measurement data Lead field matrix.

\begin{figure}[htb!]
	\centering
	\includegraphics[width=\textwidth]{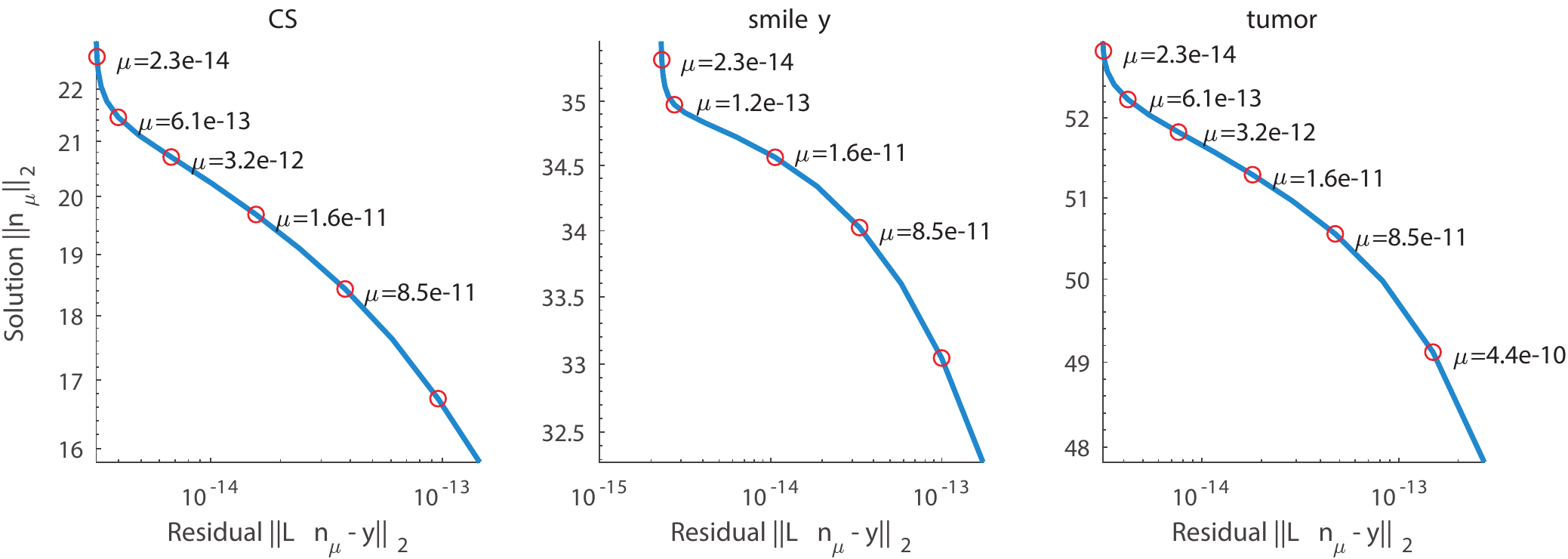}
	\caption{\textbf{Selection of the regularization parameter by the L-curve method.}\label{fig:L-Tik}
Each image shows a log-log plot of the residual functionals versus the solution norm in dependence
of the regularization parameter. The  L-curve method predicts an L-shaped curve and   advices  to select a regularization parameter close to the corner of the L-curve.
In our  results  we empirically found that taking the regularization parameter slightly  smaller than advised by  the L-curve method yield more accurate results.}
\end{figure}

For the  presented results, we use  three  different  magnetic  particle
distributions which, together with the corresponding
measurements full data, are shown in Figure~\ref{fig:phantoms}.
Any column in the measurement data (in Figure~\ref{fig:phantoms} and below)
corresponds to a single  activation pattern and contains the data of all detectors.
The phantoms are rescaled to have  maximum value 1 and the forward operator $ \Lo$ is rescaled to have
  matrix norm $\norm{\Lo}_2=1$.
 To all data we have  added additive Gaussian noise amounting to
 a signal-to-noise  (SNR) of  \SI{80}{dB}. The corresponding  reconstructions
\begin{equation}\label{eq:qtik}
    \mnp_\mu  =  (\Lo^\trans \Lo + \mu  \Io)^{-1} \Lo^\trans \bfield
\end{equation}
with standard  quadratic Tikhonov regularization (penalized  least squares)
using  regularization parameter $\mu = 10^{-12}$ are shown in the
bottom row of Figure~\ref{fig:phantoms}. Recall that $\Lo \in   \R^{\ncoils  \nsensors \times \nvoxels }$ is the full
Lead field matrix  and $\bfield \in \R^{\ncoils  \nsensors}$ the full activation data.
Reconstructions with \eqref{eq:qtik} can be seen as benchmark  for the more sophisticated CS reconstructions applied
to less data presented below. The regularization parameter has been
chosen empirically as the trade-off between stability and smoothing.
Taking a  smaller regularization parameter would not suppress the noise (amplification) well enough,
whereas a larger regularization parameter yields to overshooting. There are many  strategies
for selecting the regularization parameter analytically \cite{Han98}.  An example is the L-curve method, where the residual $\norm{\Lo \mnp_\mu - \y}_2^2$ is plotted against the solution norm  $\norm{\mnp_\mu}$ both  on a logarithmic scale. The L-curve method predicts the  graph to be
an L-shaped curve and the regularization parameter should be
taken  close to the corner of the L-shaped curve.  The L-curves  for the considered phantoms
are  plotted in Figure {\ref{fig:L-Tik}}. Indeed, the selected $\mu = 10^{-12}$  are close to the corner in each case. Strictly taken,
the L-curve method predicts a larger regularization which we however found to yield to a slight over-smoothing
and a larger reconstruction error.

The CS forward matrix $\Mo  = \Ao \Lo$ is computed  by multiplying the full
Lead field matrix with $\Ao$.  CS measurements    have  been generated in two random
ways and one   deterministic way. In the random case, $\A$ is taken either
as  Bernoulli matrix having entries $\pm 1$ appearing with equal
probability, or  a  Gaussian matrix consisting of i.i.d.
$\mathcal{N}(0,1)$-Gaussian random  variables in each entry.
The deterministic sparse sampling is performed by choosing $\nact$ equispaced coil
activations.

 \begin{figure}[htb!]
	\centering
	\includegraphics[width = 0.9\columnwidth]{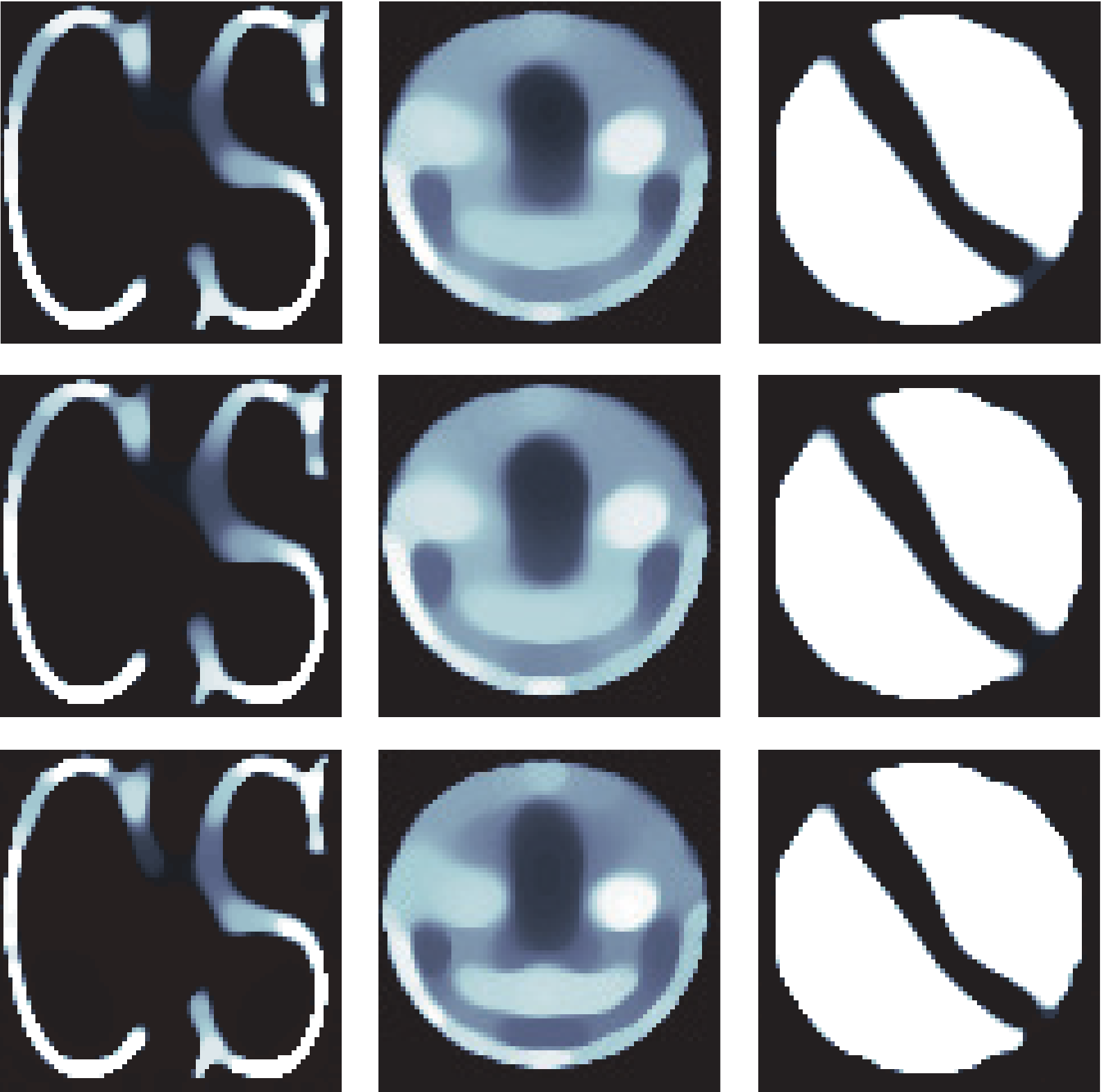}
	\caption{\textbf{Reconstruction results with  Algorithm~\ref{alg:dr}
	from 40 coil activations.}
	 Top row: Gaussian activation matrix.
	 Middle row: Bernoulli activation matrix.
	 Bottom row:  Deterministic activation pattern.}
	\label{fig:1stage}
\end{figure}

\begin{table}[htb!]
\centering
    \begin{tabular}{|c|c|c|c|}
          \toprule
           & relative RMSE $(\%)$ & SNR (dB) & correlation \\
         \toprule
        \multicolumn{4}{|c|}{ 40 Gaussian activations (proposed algorithm)}
        \\
        \midrule
         CS-phantom        &  0.75   & 2.51   & 0.92  \\
         smiley-phantom    &  0.40   & 8.05   & 0.87      \\
         tumor-phantom     &  0.16   & 15.69   & 0.98   \\
        \midrule
        \multicolumn{4}{|c|}{ 40 Bernoulli activations (proposed algorithm) }
        \\
        \midrule
         CS-phantom        &  0.78  & 2.15  & 0.88 \\
         smiley-phantom    &  0.39  &  8.11  & 0.87     \\
         tumor-phantom     &  0.17  & 15.56  & 0.98   \\
        \midrule
         \multicolumn{4}{|c|}{ 40 deterministic activations (proposed algorithm) }
        \\
        \midrule
         CS-phantom        &  0.66  & 3.59  & 0.92 \\
         smiley-phantom    &  0.35  &  9.16 & 0.90     \\
         tumor-phantom     &  0.17  & 15.38  & 0.98  \\
        \midrule
         \multicolumn{4}{|c|}{ 120 activations (Tikhonov regularization) }
        \\
         \midrule
         CS-phantom        &  1.20  & -1.59  &  0.69 \\
         smiley-phantom    &  0.47  & 6.45  & 0.80  \\
         tumor-phantom     &  0.35  & 8.99 & 0.89  \\
         \bottomrule
    \end{tabular}
    \caption{Evaluation metrics for Gaussian measurements (row 1)
    Bernoulli measurements (row 2) and  deterministic activation pattern (row 3)
    using 40  activation patterns. The bottom row shows the same
    evaluation metrics  for Tikhonov regularization for full activations.}
    \label{tab:errors}
\end{table}

\subsection{Reconstruction results}

The system matrix  $\Mo$ is rescaled to have  matrix norm 1 and we use the parameter  setting
$\mu  = 4 \cdot 10^{-13}$, $\alpha  = 10^{-14}$, $s=1$, $n_{\rm max}=1$ and
$N_{\rm iter} =50$. Note that the parameter $\mu$ has a similar role  for
each iterative step as in quadratic Tikhonov regularization {\eqref{eq:qtik}}.
Choosing it too small causes noise amplification.
Choosing it too large causes oversmoothing, which however will be reduced during the iteration.
Therefore, taking it somewhat larger than
the value found in Tikhonov case seems reasonable.  The regularization parameter
$\alpha$ has been selected empirically to yield good numerical performance.

We observed that after 50 iterations the reconstructed MNP  stagnates
which indicates that the  iterates are close to the minimizer of the Tikhonov functional.
The reconstructed MNP distributions
using   Algorithm~\ref{alg:dr}   from 40 coil activations are  shown in
Figure~\ref{fig:1stage}.
Each reconstruction  takes about   50 seconds in Matlab R 2017a  on a MacBook Pro (2016)
with 2.9 GHz Intel Core i7 processor.

To quantitatively evaluate the reconstruction results, we
compute the relative root mean squared error (RMSE)
$\norm{\mnp - \mnp_{\rm rec}}/\norm{\mnp}$,  the SNR $20 \lg (\norm{\mnp} / \norm{\mnp - \mnp_{\rm rec}})$,
and the Pearson correlation coefficient  $\operatorname{cov}( \mnp, \mnp_{\rm rec} ) / \sqrt{\operatorname{var}( \mnp ) \operatorname{var}( \mnp_{\rm rec} )}$
with $\operatorname{cov}$ denoting the covariance and $\operatorname{var}$ the variance.
The results are shown in Table~{\ref{tab:errors}}.
For all  phantoms and for any evaluation metric, the proposed algorithm  clearly  outperforms
quadratic Tikhonov regularization. More specifically, the  RSME for the proposed algorithm applied with 40 activations yields
half of the  RSME  compared to quadratic Tikhonov regularization with 120 activations.
For the CS-phantom and the smiley-phantom, the deterministic scheme
performs best and decreases the reconstruction error by more than $\SI{10}{\percent}$
compared to the random schemes. For the tumor phantom, all
sampling schemes yield comparable performance.

 \begin{figure}[htb!]
	\centering
	\includegraphics[width = \columnwidth]{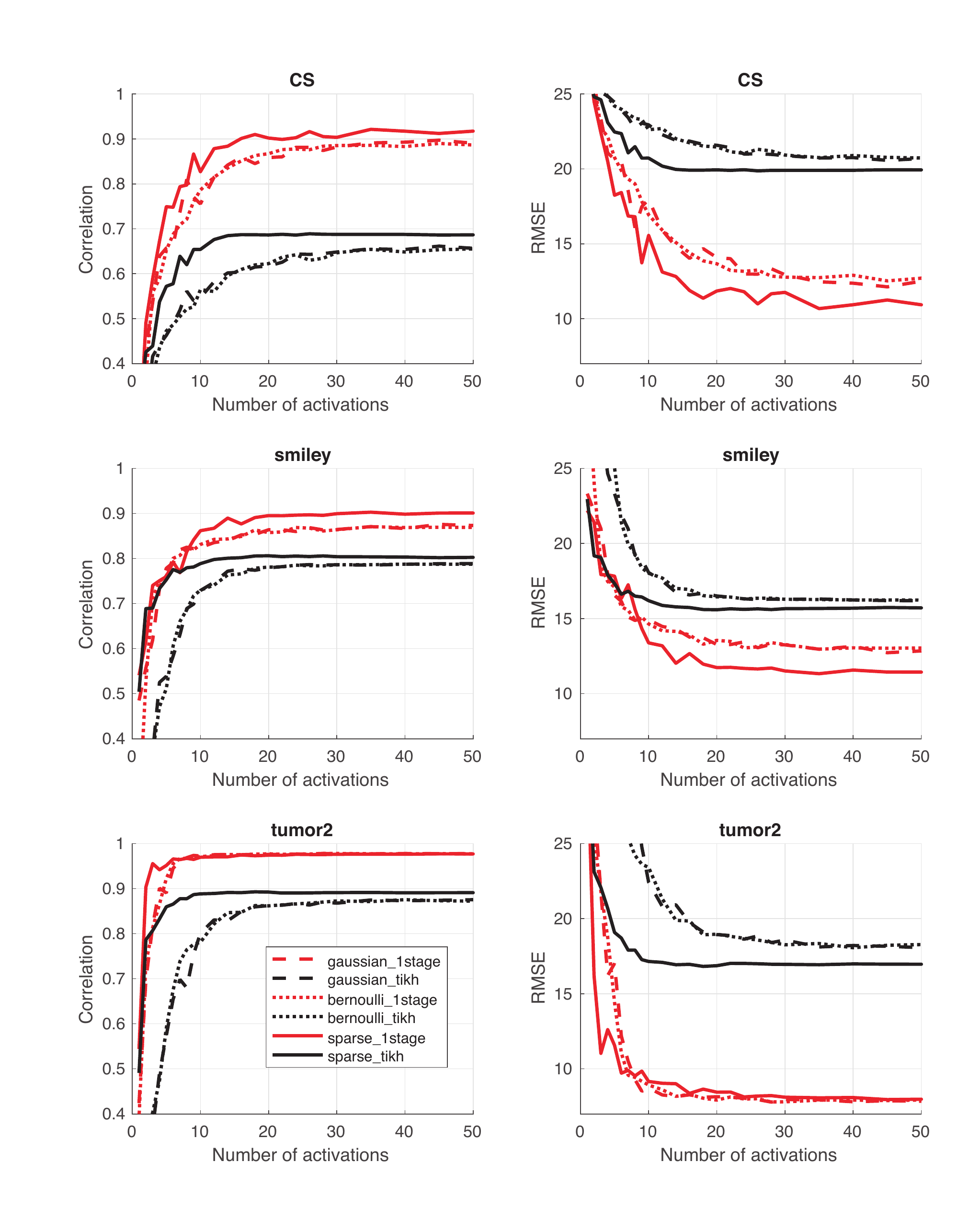}
	\caption{\textbf{Quality evaluation in dependence of the number of coil activations.}\label{fig:comparison}
Left: correlation between reconstructions and the true phantoms.
Right: RMSE of the reconstructions. }
\end{figure}

Figure~\ref{fig:comparison} shows a convergence study (for the correlation and the RMSE)
in dependence of the number of activations. We see that the Douglas-Rachford algorithm
constantly outperforms  Tikhonov regularization and significantly reduces the RSME and increases the correlation.  For the CS and the smiley phantom, the deterministic sampling pattern has a significantly smaller RMSE than the random schemes.
We address this behavior to the strong smoothing effect of the  forward  operator, which removes most high frequency components in the data. In the full forward matrix
$\A \,  \LL $, the  incoherence contained in the measurements matrix $\A$ is annihilated
by the smoothing effect of $\A$. Finding particular activation patterns that outperform the
sparse sampling scheme is a nontrivial issue and will be investigated in future work.

\section{Conclusion}
\label{sec:conclusion}

The use of multiple coil activation patterns in magnetorelaxometry imaging is time consuming and requires performing  several consecutive measurements.
It is therefore desirable to make the number of coil activations
as small as possible, while keeping high spatial resolution. For that purpose, we investigated
CS  strategies in this paper.
We compared  Gaussian random activations,
Bernoulli activations and deterministic sparse sampling schemes.
For actual image reconstruction, we applied Douglas-Rachford splitting to the sparse  Tikhonov functional; see Algorithm~\ref{alg:dr}.
For a small number of coil activations, the random
schemes slightly outperform  the deterministic scheme for the smiley-phantom and tumor
phantom. For a large number of activations, the  deterministic
sampling scheme clearly performed better  than the  random sampling patterns.
  We explain this behavior by the severe ill-conditioning  of the Lead field matrix,
  which is more severe for  a large number of activations. Only for  a small number of activations,
  where  the resolution is poor anyway, the incoherence at the low frequencies  is not destroyed
  by the smoothing effect of the Lead field matrix.
The reconstruction results clearly demonstrate the  proposed
Algorithm~\ref{alg:dr} significantly outperforms standard algorithms such as
Tikhonov regularization. Depending on the complexity of the phantom between 6 (tumor) and  30 activations (CS) are
sufficient for the used framework and  further increasing
the  number of activations only decreases the RMSE by a few percent; see Figure~\ref{fig:comparison}.

Several interesting research directions following this  work are  possible.
First, we can replace the inner
TV minimization in the single-stage approach by a different algorithm
which  should   accelerate Algorithm~\ref{alg:dr}.
Second, the derivation of theoretical error  estimates  for the single-stage approach
is of significant interest.  Results in that direction can advise which type of MRX measurements
are best  to obtain  accurate results for certain phantom classes.
Moreover, the derivation of adaptive compressed  sensing strategies
for online monitoring    is of significant interest. Advising optimal coil  activations
 given previous  activations is a practically important challenge that will benefit from
 theoretical error estimates, numerical simulations as well as real-world
 experiments.  Such issues will be addressed in future work. In this paper,
 we investigated standard random compressed sensing schemes (using Gaussian and Bernoulli
 activation patterns and sparse sampling).  Using more problem adapted and task
 oriented measurement design  we expect to  derive improved  coil activation schemes.

\section*{Acknowledgment}
This work has been supported by the Austrian Science Fund (FWF),
project P 30747-N32, and the German Science Foundation (DFG)
within the priority program COSIP, project Cos-MRXI (BA 4858/2-1).

\end{document}